\documentclass[11pt openany]{article}
\usepackage{graphics,graphicx,epsfig}
\usepackage{amsmath,amsfonts,amssymb,amstext,amsthm,amscd,mathrsfs}
\usepackage[spanish,english]{babel}
\usepackage[all]{xy}
\usepackage{multicol}

\theoremstyle{definition}

\theoremstyle{plain}
\newtheorem{theo}{Theorem}

\newcommand{\noi}{{\noindent}}
\newcommand{\su} {\mathbf{U}}
\newcommand{\gl} {\mathbf{GL}}
\newcommand{\Z}{\mathbb{Z}}
\newcommand{\C}{\mathbb{C}}
\newcommand{\R}{\mathbb{R}}
\newcommand{\CP}{\mathbb{CP}}
\newcommand{\T} {\mathbb{T}}
\newcommand{\vdu}{\mathcal{N}}
\newcommand{\vd} {\mathcal {M}}
\newcommand{\lv} {LV$-$M}

\newcommand{\Lam}{\pmb{\Lambda}}
\newcommand{\ZLam}{Z^{^\C}}

\newcommand{\s}{\mathbb{S}}
\newcommand{\disc}{\mathbb{D}}

\newcommand{\p} {\mathcal{P}}

\newcommand{\te}{\mathcal{E}}

\newcommand{\gato}{\mathop{\pmb{\#}}\limits}

\begin{document}
\title{Some Open Book and Contact Structures on Moment-Angle Manifolds.}
\author{Yadira Barreto, Santiago L\'opez de Medrano and \\ Alberto Verjovsky.\thanks{The authors were partially supported by projects PAPIIT-DGAPA IN103914 and IN108112.} }
\date{}
\maketitle
\vspace{-0.2in}
\begin{flushright}
\emph{In memory of Samuel Gitler}
\end{flushright}
\begin{abstract}
\noi I. We construct open book structures on all moment-angle ma\-ni\-folds\- and describe the topology of their leaves and bindings under certain restrictions. II. We also show, using a recent deep result about contact forms due to Borman, Eliashberg and Murphy \cite{BEM}, that every odd-dimensional moment-angle manifold admits a contact structure. This contrasts with the fact that, except for a few, well-determined cases, even-dimensional ones do not admit symplectic structures. We obtain the same results for large families of more general intersections of quadrics.
\end{abstract}
\vspace{0.2in}

\noi {\bf Key Words:} Open book decomposition, moment-angle manifolds,  contact structures.
\section*{Introduction}
The topology of generic intersections of quadrics in $\R^n$ of the form:
$$
\sum_{i=1}^n\lambda_ix^2_i=0,\quad\quad\sum_{i=1}^n x^2_i=1,
$$
where $\lambda_i\in\R^k$, $i=1,\dots,n$, has been studied for many years (\cite{SLM}, \cite{LoGli}, \cite{SLM2}, \cite{GL}). If $\Lam=\left(\lambda_1,\cdots, \lambda_n  \right)$, we will denote this variety by $Z=Z(\Lam)$.\\

\noi We will always assume the following generic condition, known as \emph{weak hyperbo\-licity} and equivalent to the smoothness of $Z$:\\

\noi\emph{If $J\subset\{1,\dots,m\}$ has $k$ or fewer elements then the origin is not in the convex hull of the $\lambda_i$ with $i\in J$.}\\

\noi A crucial feature of these manifolds is that they admit natural group actions: all of them admit $\Z^{^n}_2$ actions by changing the signs of the coordinates.\\ 

\noi Their complex versions in $\C^n$, which we denote by $Z^{\C}$, 

$$
\sum_{i=1}^n\lambda_i|z_i|^2=0,\quad\quad \sum_{i=1}^n|z_i|^2=1,
$$
\noi (now known as \emph{moment-angle manifolds}) admit natural actions of the $n$-torus $\T^n$. The quotient can be identified in both cases with the polytope $\p$ given by

$$
\sum_{i=1}^n\lambda_ir_i=0,\quad\quad\sum_{i=1}^nr_i=1,\quad\quad r_i\geqslant 0.
$$
that determines completely the varieties (so we can use the notations $Z(\p)$ and $Z^{^\C}(\p)$ for them) as well as the actions. The weak hyperbolicity condition implies that $\p$ is a simple polytope and any simple polytope  can be realized as such a quotient. The facets of $\p$ are its non-empty intersections with the coordinate hyperplanes. If all such intersections are non-empty $Z$ and $Z^{^\C}$ fall under the general concept of \emph{generalized moment-angle complexes} (\cite{BBCG}).\\

\noi If we take the quotient of $\ZLam$ by the scalar action of $\s^{^1}$: 
$$
\vdu(\Lam)=\ZLam/\s^{^1},
$$

\noi we obtain a compact, smooth manifold  $\vdu(\Lam)\subset \CP^{^{n-1}}$.\\

\noi When $k$ is even, $\vdu(\Lam)$ and $\ZLam\times\s^{^1}$ have natural complex structures and so does $\ZLam$ itself when $k$ is odd, but admit symplectic structures only in a few well-known cases (\cite{LV}, \cite{Meer}). Those complex manifolds called now $\lv-m$-\textit{manifolds} .\\

\noi An open book construction was used to describe the topology of $Z$ for $k=2$ in some cases not covered by Theorem 2 in \cite{SLM} (see remark on page 281 and \cite{GL}). In \cite{LoGli} it is a principal technique for studying the case $k>2$. In section I-1 we recall this construction, underlining the case of moment-angle manifolds:\\

\noi \textit{If $\p$ is a simple convex polytope and $F$ one of its facets, $Z^{\C}(\p)$ admits an open book decomposition with binding $Z^{^\C}(F)$ and  trivial mo\-no\-dro\-my.}\\

\noi When $k=2$, the varieties $Z$ and $\ZLam$ can be put in a normal form given by an \emph{odd cyclic partition} (see section I-1) and they are  diffeomorphic to a triple product of spheres or to the connected sum of sphere products (see \cite{SLM}, \cite {GL} and section I-4). Using the same normal form, we give a topological description of the leaves of their open book decompositions  (improving the statement in \cite{BLV} and \cite{SLM2}) which is complete in the case of moment-angle manifolds:\\
\newpage
\noi \textit{The leaf of the open book decomposition of $\ZLam$ is the interior of:}
\begin{itemize}
\item [a)]\noi  \textit{a product $\s^{2n_2-1}\times\s^{2n_3-1}\times\disc^{2n_1-2}$,}

\item [b)]\noi  \textit{a connected sum along the boundary of products of the form $\s^p\times\disc^{2n-p-4}$,}

\item [c)]\noi  \textit{in some cases, there may appear summands of the form:}
\vspace{0.025in}

\noi \textit{a punctured product of spheres $\s^{2p-1}\times\s^{2n-2p-3}\backslash\disc^{2n-4}$ or}
\vspace{0.025in}

\noi \textit{the exterior of an embedding $\,\,\,\,\,\,\,\s^{2q-1}\times\s^{2r-1} \subset \,\,\,\s^{2n-4}.$}
\end{itemize}

\noi The precise result (Theorem 2 in section I-1) follows from Theorem 3 in section I-4, a general theorem that gives the topological description of the \emph{half} real varieties $Z_{_+}=Z\cap\{x_1\geqslant 0\}$, complementing \cite{SLM}, and requires additional dimensional and connectivity hypotheses that should be avoidable using the methods of \cite{GL}. Some of the proofs follow directly from the result in \cite{SLM}, while other ones require the use of some parts of its proof. All these manifolds with boundary are also generalized moment-angle complexes.\\

\noi In part II, using a recent deep result about contact forms due to Borman, Eliashberg and Murphy \cite{BEM}, we show that every odd dimensional moment-angle manifold admits a contact structure. This is surprising since even dimensional moment-angle manifolds admit symplectic structures only in a few well-known cases. We also show this for large families of more general odd-dimensional intersections of quadrics by a different argument.\\

\noi The original aim of this work (see \cite{BLV}) was to construct contact structures on odd-dimensional moment angle-manifolds based on the open book structures, using results of E. Giroux  and  E. Giroux and J-P. Mohsen, that relate these two structures (\cite{Giroux}, \cite{GiMo}). However, it turned out that there does not exist a contact form which is supported in the open book decompositions, as in Giroux's theorem, because the pages are not Weinstein manifolds (\cite{Weinstein}). What we do show is that all moment-angle mani\-folds admit both structures separately. We think that these results are interesting in themselves. In \cite{BV} the first and  third authors give a different construction, in some sense more explicit, of contact structures, not on moment-angle manifolds but on certain non-diagonal generalizations of moment-angle manifolds of the type that has been studied by  G\'omez Guti\'errez and the second author in \cite{GL}. It consists in the construction of a positive confoliation which is conductive and then uses the heat flow method described in \cite{AltWu}.\\

\noi Samuel Gitler taught us how to look at our varieties from the more general pers\-pective of generalized moment-angle complexes. This opened the way to the solution of many old and new problems, starting with the article \cite{LoGli} fo\-llowed by many others, including \cite{GL},  \cite{SLM2}, the present one and surely many more to come. We miss him very much personally and we miss very much the lessons he did not have time to teach us.\\

\section*{Part I. Open book structures.}

\subsection*{I-1. Construction of the open books.}

\noi Let $\Lam'$ be obtained from $\Lam$ by adding an extra $\lambda_1$ which we interpret as the coefficient of a new extra variable $x_0$, so we get the variety $Z'$:
$$
\lambda_1\left(x^2_0+x^2_1\right)+\sum_{i>1}\lambda_ix^2_i=0,\,\,\,\,\,\,\,\,\,\,(x^2_0+x^2_1)+\sum_{i>1}x^2_i=1.
$$

\noi Let $Z_+$ be the intersection of $Z$ with the half space $x_1\geqslant 0$. $Z_+$ admits an action of $\Z^{^{n-1}}_2$ with quotient the same $\p$: $Z_+$ can be obtained by reflecting $\p$ on all the coordinate hyperplanes except $x_1=0$. $Z_+$ is a manifold with boundary $Z_0$ which is the intersection of $Z$ with the subspace $x_1=0$.\\

\noi Consider the action of $\s^1$ on $Z'$ by rotation of the coordinates  $\left(x_0,x_1\right)$ . This action fixes the points of $Z_0$ and all its other orbits cut $Z_+$ transversely in exactly one point. So $Z'$ is the open book with binding $Z_0$, page $Z_+$ and trivial monodromy:

\begin{theo}\label{bookR}
i) Every manifold $Z'$ is an open book with trivial monodromy, bin\-ding $Z_0$ and page $Z_+$.\\

\noi ii) If $\p$ is a simple convex polytope and $F$ is one of its facets, there is an open book decomposition of $Z^{^\C}(\p)$ with binding $Z^{^\C}(F)$ and trivial mono\-dromy.
\end{theo}

\noi (ii) follows because if we write the equations of $Z^{^\C}(\p)$ in real coordinates, we get terms $\lambda_i (x_i^2+y_i^2)$ so each $\lambda_i$ appears twice as a coefficient and $Z^{^\C}(\p)$ is a variety of the type $Z'$ in several ways. It is then an open book with binding the manifold $Z^{^\C}_0(\p)$ obtained by taking $z_i=0$.\\

\noi When $k=2$ it can be assumed $\Lam$ is one of the following normal forms (see \cite{SLM}): Take $n=n_1+\dots+n_{2\ell+1}$ a partition of $n$ into an odd number of positive integers. Consider the configuration $\Lam$ consisting of the vertices of a regular polygon with $(2\ell+1)$ vertices, where the $i$-th vertex in the cyclic order appears with multiplicity $n_{i}$. \\
\begin{center}
\begin{figure}[h]
\vspace{-0.3in}\centerline{\includegraphics[height=3.5cm]{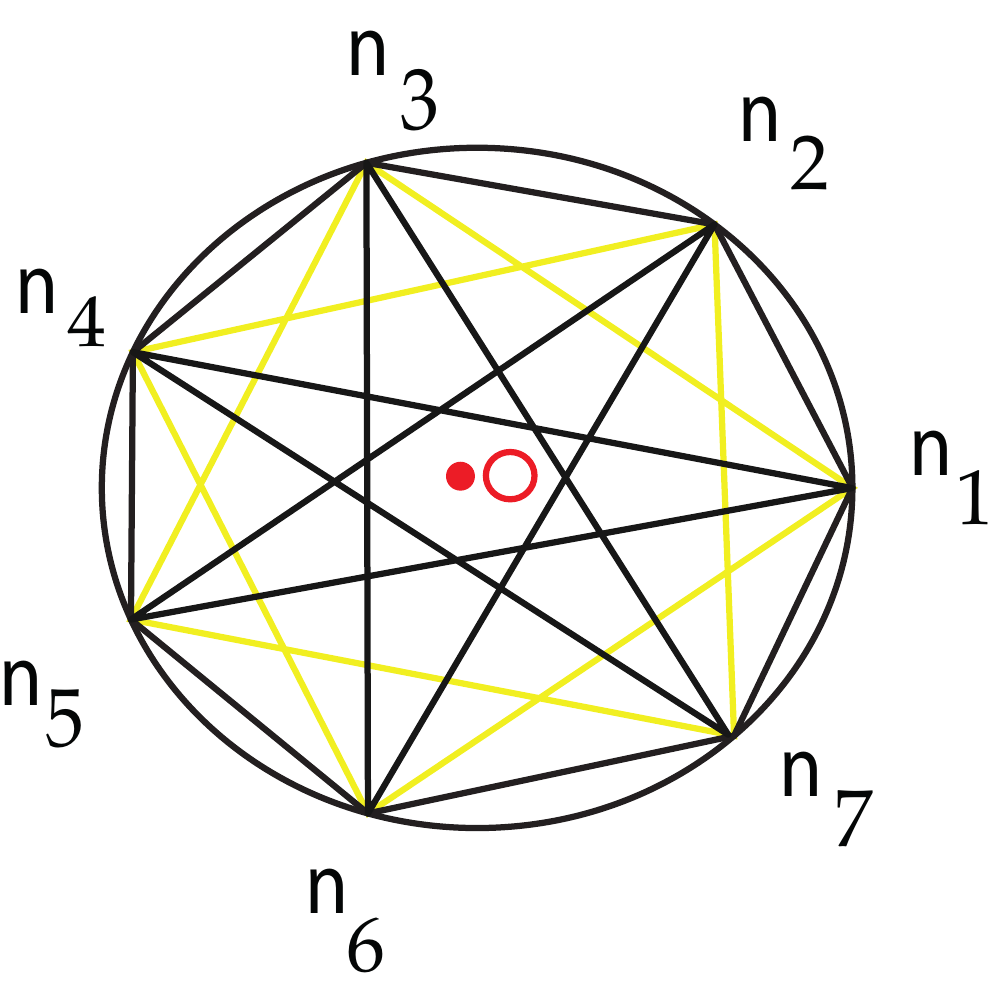}}
\end{figure}
\end{center}   
\vspace{-0.35in}
\noi The topology of $Z$ and $\ZLam$ can be completely described in terms of the numbers $d_i=n_i+\dots+n_{i+\ell-1}$, i.e., the sums of $\ell$ consecutive $n_i$ in the cyclic order of the partition (see \cite{SLM}, \cite{GL} and section I-1):
\begin{center}
\begin{figure}[h]
\centerline{\includegraphics[height=3.5cm]{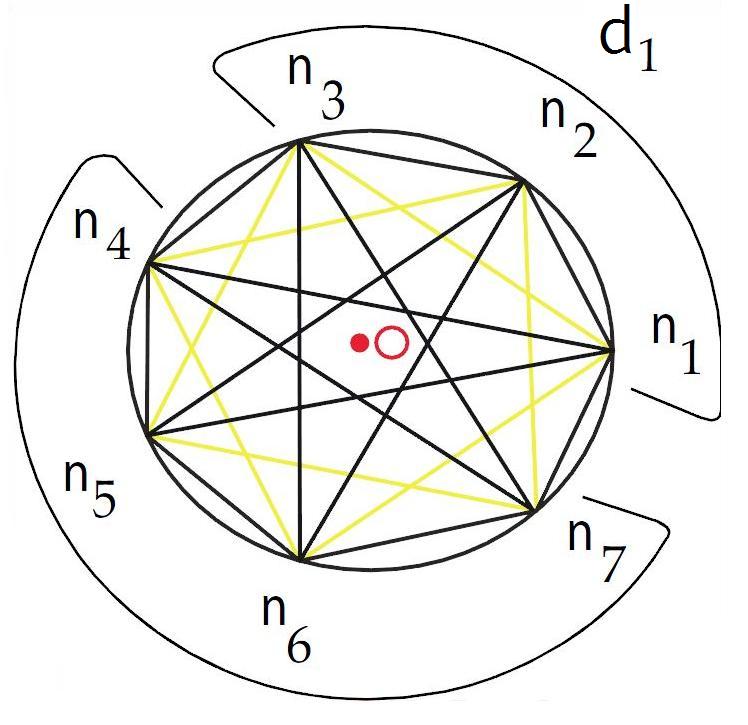}}
\end{figure}
\end{center}  
\vspace{-0.3in}
\noi For $\ell=1$: $Z=\s^{n_1-1}\times\s^{n_2-1}\times\s^{n_3-1},\,\,\,Z^\C=\s^{2n_1-1}\times\s^{2n_2-1}\times\s^{2n_3-1}$.\\

\noi For $\ell>1$: $Z=\gato_{j=1}^{2\ell+1}\left(\s^{d_i-1}\times\s^{n-d_i-2}\right),\,\,\, Z^\C=\gato_{j=1}^{2\ell+1}\left(\s^{2d_i-1}\times\s^{2n-2d_i-2}\right)$.\\

\noi Now we have a similar description of the topology of the leaves in all moment-angle manifolds, where $\coprod$ denotes connected sum along the boundary and $\te_{2n_2-1,2n_4-1}^{2n-4}$ is the exterior of $\s^{2n_2-1}\times\s^{2n_5-1}$ in $\s^{2n-4}$ (see section I-3):

\begin{theo}\label{pageC}
Let $k=2$, and consider the manifold $Z^{^\C}$ corresponding to the cyclic partition $n=n_1+\dots+n_{2\ell+1}$. Consider the open book decomposition of $Z^{^\C}$ corresponding to the binding at $z_1=0$, as given by Theorem 1. Then the leaf of this decomposition is diffeomorphic to the interior of:

\begin{itemize}
\item [a)] If $\ell=1$, the product 
$$
\s^{2n_2-1}\times\s^{2n_3-1}\times\disc^{2n_1-2}.
$$

\item [b)] If $\ell>1$ and $n_1>1$, the connected sum along the boundary of $2\ell+1$ manifolds:
$$
\coprod_{i=2}^{\ell+2}\left(\s^{2d_i-1}\times\disc^{2n-2d_i-3}\right)\coprod\coprod_{i=\ell+3}^1\left(\disc^{2d_i-2}\times\s^{2n-2d_i-2}\right).
$$

\item [c)] If $n_1=1$ and $\ell>2$, the connected sum along the boundary of $2\ell$ manifolds:
$$
\coprod_{i=3}^{\ell+1}\left(\s^{2d_i-1}\times\disc^{2n-2d_i-3}\right)\coprod\coprod_{i=\ell+3}^1\left(\disc^{2d_i-2}\times\s^{2n-2d_i-2}\right)
$$
$$
\coprod\left(\s^{2d_2-1}\times\s^{2d_{\ell+2}-1}\backslash\disc^{2n-4}\right).
$$

\item [d)] If $n_1=1$ and $\ell=2$, the connected sum along the boundary of two manifolds:
$$
\left(\s^{2d_2-1}\times\s^{2d_4-1}\backslash \disc^{2n-4}\right)\coprod\te_{2n_2-1,2n_5-1}^{2n-4}.
$$

\end{itemize}
\end{theo} 

\noi Theorem 2 will follow from its real version (see Theorem 3). It follows also that in cases c) and d) the product of the leaf with an open interval is diffeomorphic to the interior of a connected sum along the boundary of the type of case b).\\

\noi For $k>2$, the topology of moment-angle manifolds and their leaves is much more complicated and it seems hopeless to give a complete description of them: they may have non-trivial triple Massey products (\cite{Ba}), any amount of torsion in their homology (\cite{BM}) or may be a different kind of open books (\cite{LoGli}). Nevertheless, it is plausible that a description of their leaves as above may be possible for large families of them in the spirit of \cite{LoGli}.\\

\noi The manifold $\vdu(\Lam)$, defined in the introduction, also admits an open book decomposition, since the $\s^1$ action on the first coordinate commutes with the diagonal one.\\

\noi Let
$$
\pi_{_{\Lambda}}:\ZLam\to\vdu(\Lam),
$$
denote the canonical projection.\\

\noi Consider now the open book decomposition of $Z^{^\C}$ described above, co\-rres\-pon\-ding to the variable $z_{_1}$. If $\Lam_{_0}$ is obtained from $\Lam$ by removing $\lambda_{_1}$ it is clear that the diagonal $\s^{^1}$-action on $Z^{^\C}$ has the property that each orbit intersects each page in a unique point and at all of its points this page is intersected tranversally by the orbits. This implies that the restriction of the canonical projection $\pi_{_{\Lambda}}$ to each page is a diffeomorphism onto its image $\vdu(\Lam)-\vdu(\Lam_{_0})$.\\

\noi For $k$ even we therefore obtain, since $\vdu(\Lam)-\vdu(\Lam_{_0})$ has a complex structure:\\

\noi \emph{For $k$ even, the page of the open book decomposition of $Z^{^\C}(\pmb{\Lambda})$ in Theorem 2 with binding $Z^{^\C}_{_0}(\pmb{\Lambda}_{_0})$ admits a natural complex structure which makes it biholomorphic to $\vdu(\pmb{\Lambda})-\vdu(\pmb{\Lambda}_{_0})$.}\\

\noi \emph{For $k$ odd, both the whole manifold and the binding admit natural complex structures.}\\

\noi So we have a very nice open book decomposition of every moment-angle ma\-nifold. Unfortunately, it does not have the right properties to help in the cons\-truction of contact structures on them.\\

\noi The topology of these manifolds and of the leaves of their foliations is more complicated, even for $k=2$, and only some cases have been described (see \cite{LV} for the simpler ones).

\subsection*{I-2. Homology of intersections of quadrics and their halves.}\label{I 2}

We recall here the results of \cite{SLM}, whose proofs are equally valid for any intersection of quadrics and not only  for $k=2$.\\

\noi Let $Z=Z(\Lam)\subset \R^n$ as before, $\p$ its associated polytope  and $F_1,\dots,F_n$ the intersections of $\p$ with the coordinate hyperplanes $x_i=0$ (some of which might be empty).\\

\noi Let $g_i$ be the reflection on the $i$-th coordinate hyperplane and for $J\subset\{1,\dots,n\}$ let $g_J$ be the composition of the $g_i$ with $i\in J$. Let also $F_J$ be the intersection of the $F_i$ for $i\in  J$.\\

\noi The polytope $\p$, all its faces (the non-empty $F_J$) and all their combined reflections on the coordinate hyperplanes form a cell decomposition of $Z$. Then the elements $g_J(F_L)$ with non-empty $F_L$ generate the chain groups $C_\ast(Z)$, where to avoid repetitions one has to ask $J\cap L=\emptyset$ (since $g_i$ acts trivially on $F_i$).\\

\noi A more useful basis is given as follows: let $h_i=1-g_i$ and $h_J$ be  the product of the $h_i$ with $i\in J$. The elements $h_J(F_L)$ with $J\cap L=\emptyset$ are also a basis, with the advantage that $h_JC_\ast(Z)$ is a chain subcomplex of $C_\ast(Z)$ for every $J$ and, since $h_i$ annihilates $F_i$ and all its subfaces, this subcomplex can be identified with the chain complex $C_\ast(\p,\p_J)$, where $\p_J$ is the union of all the $F_i$ with $i\in J$. It follows that
$$
H_\ast(Z)\approx\oplus_JH_\ast(\p,\p_J).
$$

\noi For the manifold $Z_+$ we start also with the faces of $\p$, but we cannot reflect them in the subspace $x_1=0$. This means we miss the classes $h_J(F_L)$ where $1\in J$ and we get\footnote{The retraction $Z\to Z_{_+}$ induces an epimorphism in homology and fundamental group.}
$$H_\ast(Z_{_+})\approx\oplus_{1\notin J} H_\ast(\p,\p_J).$$

\noi To compute the homology of $\ZLam$ one can just take that of its real version (with each $\lambda_i$ duplicated) or directly using instead of the basis $h_J(F_L)$ with $J\cap L=\emptyset$ the basis of (singular) cells $F_L \times T_J$ (with $J\cap L=\emptyset$) where $T_J$ is the subtorus of $T^n$ which is the product of its $i-$th factors with $i\in J$. This gives the splitting 

$$
H_i(\ZLam)\approx\oplus_J H_{i-\vert J \vert}(\p,\p_J).
$$

\noi This splitting was actually proved before the previous one, but was not needed (and thus not written) in \cite{SLM}. It was rediscovered in \cite{BM}.\\

\noi These splittings have the same summands as the ones in \cite{BBCG} derived from the homotopy splitting of $\Sigma Z$. Even if it is not clear that they are the \emph{same} splitting, having two such with different geometric interpretations is most valuable.

\subsection*{I-3. The space $\te_{p,q}^m$.}
\bigskip

\noi Consider the standard embedding of $\s^{p}\times \s^q$ in $\s^m$, $m>p+q$ given by 

$$\s^{p}\times \s^{q}\subset \R^{p+1}\times \R^{q+1}=\R^{p+q+2}\subset \R^{m+1}.$$
whose image lies in the $m$-sphere of radius $\sqrt{2}$.\\

\noi We will denote by $\te_{p,q}^m$ the exterior of this embedding, i.e., the complement in $S^m$ of the open tubular neighborhood $U=int(\s^{p}\times \s^q\times \disc^{m-p-q})\subset \s^{m}$. Observe that the boundary of $\te_{p,q}^m$ is $\s^{p}\times \s^q\times \s^{m-p-q-1}$ and that the classes $[\s^{m-p-q-1}]$, $[\s^{p}\times \s^{m-p-q-1}]$ and $[\s^{q}\times \s^{m-p-q-1}]$ are the ones bellow the top dimension that go to zero in the homology of $U$. By Alexander duality, the images of these classes freely generate the homology of $\te_{p,q}^m$.\\

\noi Theorem A2.2 of \cite{LoGli} tells that, under adequate hypotheses (and probably always) $\te_{p,q}^m \times \disc^1$ is diffeomorphic to a connected sum along the boundary of products of the type $\s^a \times \disc^{m+1-a}$.\\

\noi Under some conditions (and probably always), $\te_{p,q}^m$ is characterized by its boundary and its homology properties: Let $X$ be a smooth compact manifold with boundary $\s^{p}\times \s^q\times \s^{m-p-q-1}$ and $\iota$ the inclusion $\partial X\subset X$.\\

\noi \textbf{Lemma.} \textit{Assume that}
\vspace{0.1in}

(i) \textit{$X$ and $\partial X$ are simply connected.}
\vspace{0.1in}

(ii) \textit{the classes $\iota_*[\s^{m-p-q-1}]$, $\iota_*[\s^{p}\times \s^{m-p-q-1}]$ and $\iota_*[\s^{q}\times \s^{m-p-q-1}]$ freely generate the homology of $X$.}
\vspace{0.1in}

\textit{Then $X$ is diffeomorphic to $\te_{p,q}^m$.}\\

\noi \textbf{Proof:} Observe that condition (i) implies that $p,q,m-p-q-1 \ge 2$ so $dim(X)=m \ge 7$. Consider the following subset of $\partial X$:
$$K=\s^{p}\times * \times \s^{m-p-q-1} \cup * \times \s^{q}\times \s^{m-p-q-1}$$
and embed $K$ into the interior of $X$ as $K\times \{1/2\}$ with respect to a collar neighborhood $\partial X \times [0,1)$ of $\partial X$. Finally, let $V$ be a smooth regular neighborhood (\cite{H}) of $K\times \{1/2\}$ in $\partial X \times [0,1)$.\\

\noi Now, the inclusion $V \subset X$ induces an isomorphism in homology. Since the codimension of $K$ in $X$ is equal to $1 + min(p,q)\ge 3$, $X\setminus int(V)$ is simply connected and therefore an h-cobordism, so $X$ is diffeomorphic to $V$.\\

\noi Since $\te_{p,q}^m$ verifies the same properties as $X$, the above construction \textit{with the same $V$} shows that $\te_{p,q}^m$ is also diffeomorphic to $V$ and the Lemma is proved.

\subsection*{I-4 Topology of $Z$ and $Z_+$ when $k=2$}\label{I 4}

{\noi For $k=2$ and $\ell=1$ a simple computation shows that
$$
Z_+=\disc^{n_1-1}\times\s^{n_2-1}\times\s^{n_3-1}.
$$}

\noi For the case $\ell>1$ we recall here the main steps in the proof of the result about the topology of $Z$ in \cite{SLM}, underlining those that are needed to determine the topology of $Z_+$. For the cyclic partition $n=n_1+\dots+n_{2\ell+1}$ we will denote by $J_i$ the set of indices corresponding to the $n_i$ copies of the $i$-th vertex of the polygon in the normal form (see I-1). Let also $D_i=J_i\cup\dots\cup J_{i+\ell-1}$ and $\tilde{D}_i$ its complement.
\begin{center}
\begin{figure}[h]
\vspace{-0.25in}\centerline{\includegraphics[height=3.5cm]{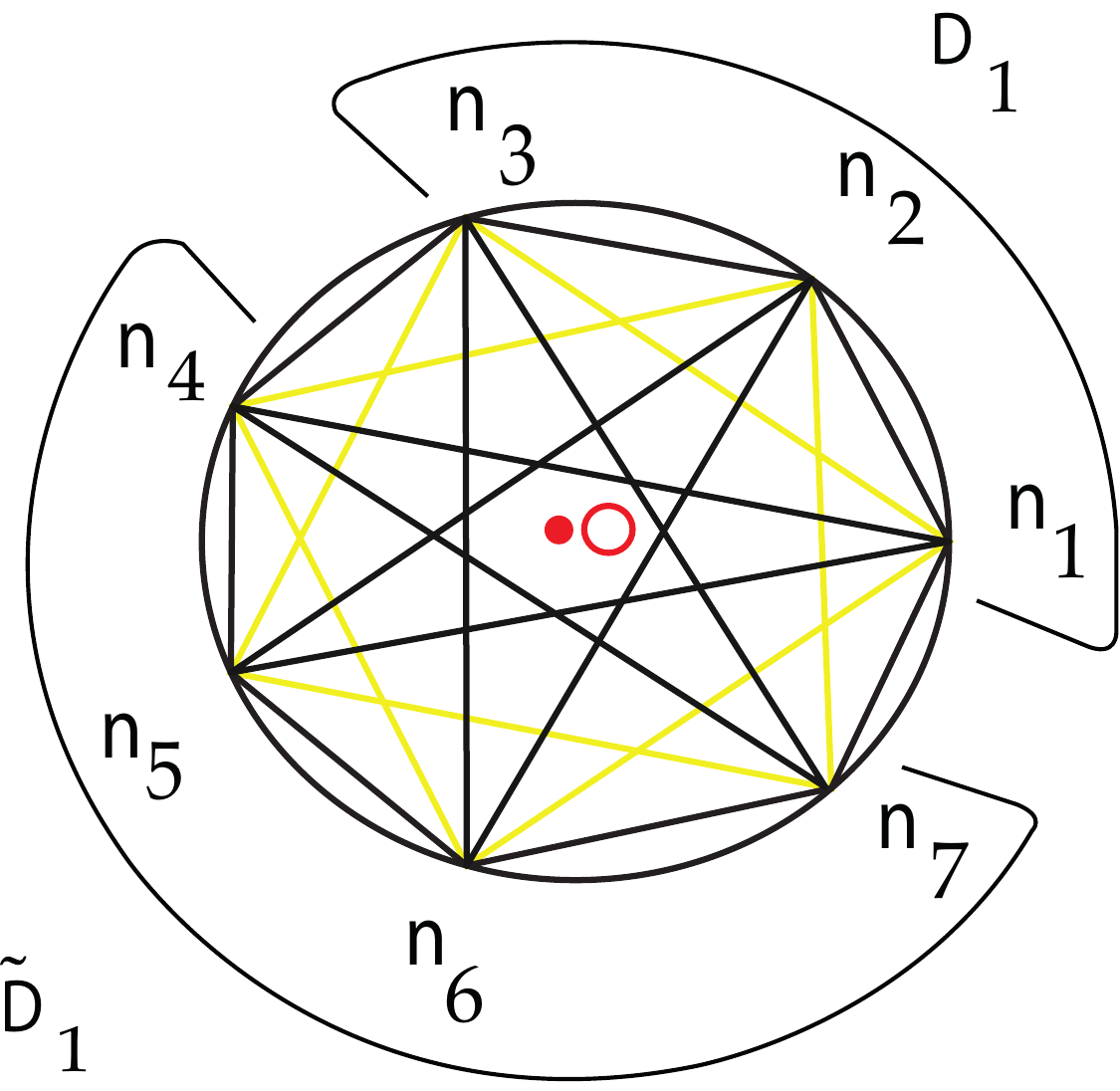}}
\end{figure}
\end{center}   
\vspace{-0.4in}
\noi It is shown in \cite{SLM} that for $k=2$, the pairs $(\p,\p_J)$ with non-trivial homology are those where $J$ consists of $\ell$ or $\ell+1$ consecutive classes, that is, those where $J$ is some $D_i$ or $\tilde{D}_i$. In those cases there is just one dimension where the homology is non-trivial and it is infinite cyclic.\\

\noi  In the case of $D_i$ that homology group is in dimension $d_i-1$ where $d_i=n_i+\dots+n_{i+\ell-1}$ is the length of $D_i$. A generator is given by the face $F_{L_i}$ where
$$L_i=\tilde{D}_i\backslash\left(\{j_{i-1}\}\cup\{j_{i+\ell}\}\right)$$
and $j_{i-1}\in J_{i-1}$, $j_{i+\ell}\in J_{i+\ell}$ are any two indices in the extreme classes of $\tilde{D}_i$ (in other words, those conti\-guous to $D_i$).\\
 
\noi $F_{L_i}$ is non empty of dimension $d_i-1$. It is not in $\p_{D_i}$, but its boundary is. Therefore it represents a homology class in $H_{d_i-1}(\p,\p_{D_i})$, which defines a gene\-rator $h_{D_i}F_{L_i}$ of $H_{d_i-1}(Z)$. Since $F_{L_i}$ has exactly $d_i$ facets it is a $(d_i-1)$-simplex so when reflected in all the coordinate subspaces containing those facets we obtain a sphere, which clearly represents $h_{D_i}F_{L_i}\in H_{d_i-1}(Z)$.\\
 
 \noi The class corresponding to $\tilde{D}_i$ is in dimension $n-d_i-2$ and is represented by the face $F_{\tilde{L}_i}$, where $\tilde{L}_i=D_i\backslash\{j\}$ and $j$ is any index in one of the extreme classes of $D_i$. It represents a generator of $H_{n-d_i-2}(Z)$, but now it is a product of spheres. For $\ell=1$ this cannot be avoided, but for $\ell>1$, with a good choice of $j$ and a surgery, it can be represented by a sphere (this also follows from \cite {LoGli}). We will not make use of this class in what follows.\\
 
 \noi The final result is that, if $\ell>1$, all the homology of $Z$ below the top dimension can be represented by embedded spheres with trivial normal bundle.\\
 
\noi Let $Z_+'$ be the manifold with boundary obtained by setting $x_0\ge0$ in $Z'$ (as defined in section I-1). Then $Z_+'$ can be deformed down to $Z_+$ by folding gradually the half-plane $x_0\ge 0, x_1$ onto the ray $x_1\ge 0$. This shows that the inclusion $Z \subset Z'_+$ induces an epimorphism in homology so one can represent all the classes in a basis of $H_\ast(Z_{_+}')$ by embedded spheres with trivial normal bundle. Those spheres can be assumed to be disjoint since they all come from the boundary $Z$ and can be placed at different levels of a collar neighborhood. Finally, one forms a manifold $Q$ with boundary by joining disjoint tubular neighborhoods of those spheres by a minimal set of tubes and then the inclusion $Q\subset Z'_+$ induces an isomorphism in homology. If $Z$ is simply connected and of dimension at least 5, then $Z'_+$ minus the interior of $Q$ is an $h$-cobordism and therefore $Z$ is diffeomorphic to the boundary of $Q$ which is a connected sum of spheres products. Knowing its homology we can tell the dimensions of those spheres:\\

\noi \emph{If $\ell>1$ and $Z$ is simply connected of dimension at least $5$, then:}
$$
Z=\gato_{j=1}^{2\ell+1}\left(\s^{d_j-1}\times\s^{n-d_j-2}\right).
$$
For the moment-angle manifold $Z^\C$ this formula gives, without any restrictions
$$
Z^\C=\gato_{j=1}^{2\ell+1}\left(\s^{2d_j-1}\times\s^{2n-2d_j-2}\right).
$$
\noi (In \cite{GL} this has recently been proved without any restrictions also on $Z$).\\

\noi The topology of $Z_+'$ is implicit in the above proof: $Z_+'$ is diffeomorphic to $Q$ and therefore it is a connected sum along the boundary of manifolds of the form $\s^p\times\disc^{n-3-p}$. Since any $Z$  with $n_1>1$ is such a $Z'$ we have:\\

\noi \emph{If $Z_0$ is simply connected of dimension at least $5$, and $\ell>1$, $n_1>1$ then:}
$$
Z_+=\coprod_{i=2}^{\ell+2}\left(\s^{d_i-1}\times\disc^{n-d_i-2}\right)\coprod\coprod_{i=\ell+3}^1\left(\disc^{d_i-1}\times\s^{n-d_i-2}\right).
$$
The clases $D_i$ and $\tilde{D}_i$ that now give no homology are the ones that contain $n_{1}$.\\

\noi The case $n_1=1$ is different. When $n_1>1$ the inclusion $Z_0\subset Z_+$ induces an epimorphism in homology (since it is of the type $Z\subset Z'_+$). This is not the case for $n_1=1$: for the partition $5=1+1+1+1+1$, the polytope $\p$ is a pentagon and an Euler characteristic computation (from a cell decomposition formed by $\p$ and its reflections) shows that $Z$ is the surface of genus $5$. Now $Z_0$ has partition $4=1+2+1$ and consists of four copies of $\s^1$. From this, $Z_+$ must be a torus minus four disks that can be seen as the connected sum of a sphere minus four disks (all whose homology comes from the boundary) and a torus that carries the homology not coming from the boundary.\\

\noi In general, when $n_1=1$ $Z_0$ is given by a normal form with $2\ell-1$ classes, has $4\ell-2$ homology generators below the top dimension, only half of which survive in $Z_+$. But $Z_+$ has $2\ell+1$ homology generators, so two of them do not come from its boundary and actually form a handle. To be more precise, the removal of the element $1\in J_1$ allows the opposite classes $J_{\ell+1}$ and $J_{\ell+2}$ to be joined into one without breaking the weak hyperbolicity condition.
\begin{center}
\begin{figure}[h]
\vspace{-0.25in}\centerline{\includegraphics[height=3.5cm]{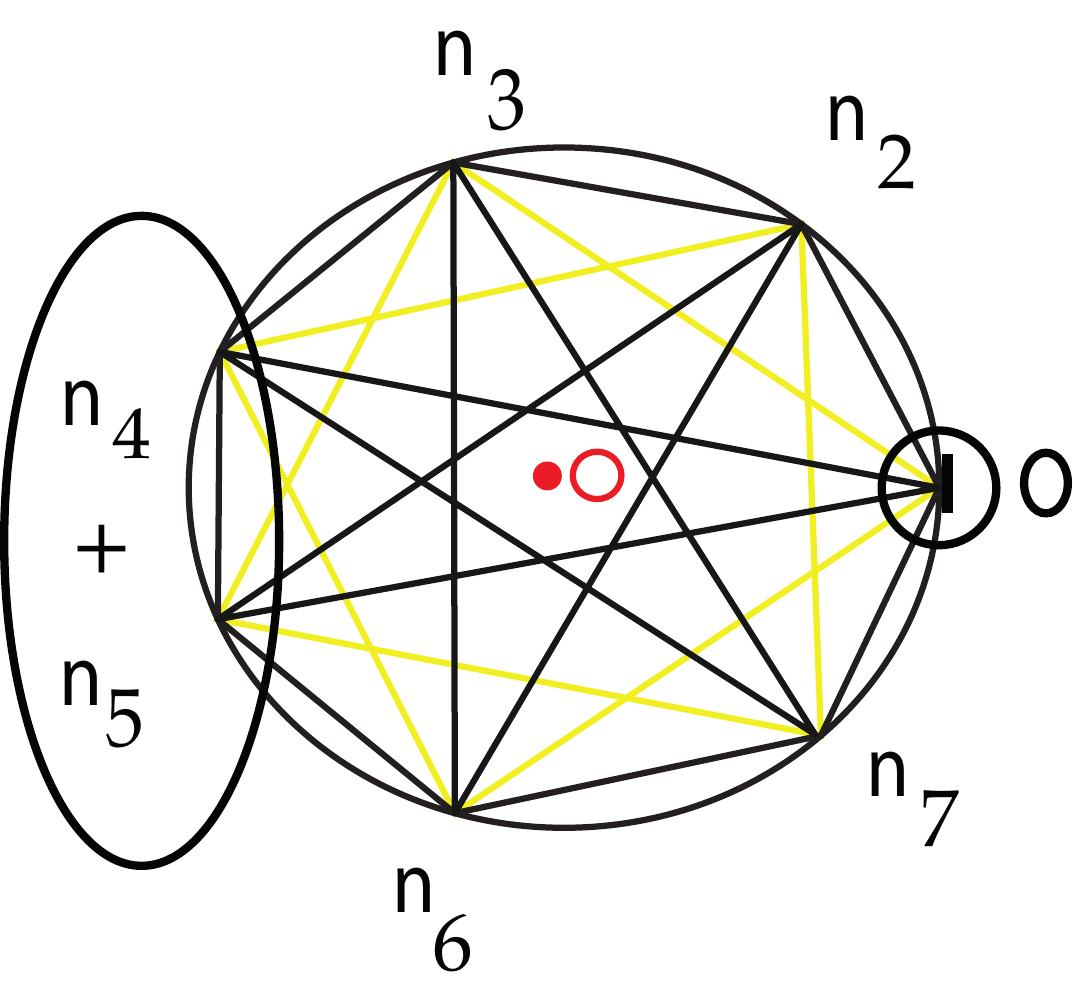}}
\end{figure}
\end{center}   
\vspace{-0.35in}
\noi Therefore $Z_0$ has fewer such classes and $D_2=J_2\cup\dots\cup J_{\ell+1}$, which gives a generator of $H_\ast(Z_{_+})$, does not give anything in $H_\ast(Z_{_0})$ because there \emph{it is not a union of classes}  (it lacks the elements of $J_{\ell+2}$ to be so).\\

\noi The two classes in $H_\ast(Z_{_+})$ missing in $H_\ast(Z_{_0})$ are thus those corresponding to $J=D_2$ and $J=D_{\ell+2}$; all the others contain both $J_{\ell+1}$ and $J_{\ell+2}$ and thus live in $H_\ast(Z_{_0})$.\\

\noi As shown above, these two classes are represented by embedded spheres in $Z_+$ with trivial normal bundle built from the cells $F_{L_2}$ and $F_{L_{\ell+2}}$ by reflection. Now $F_{L_2}\cap F_{L_{\ell+2}}$ is a single vertex $v$, all coordinates except $x_1$, $x_{j_{\ell+1}}$, $x_{j_{\ell+2}}$  being $0$.
\begin{center}
\begin{figure}[h]
\vspace{-0.25in}\centerline{\includegraphics[height=3.5cm]{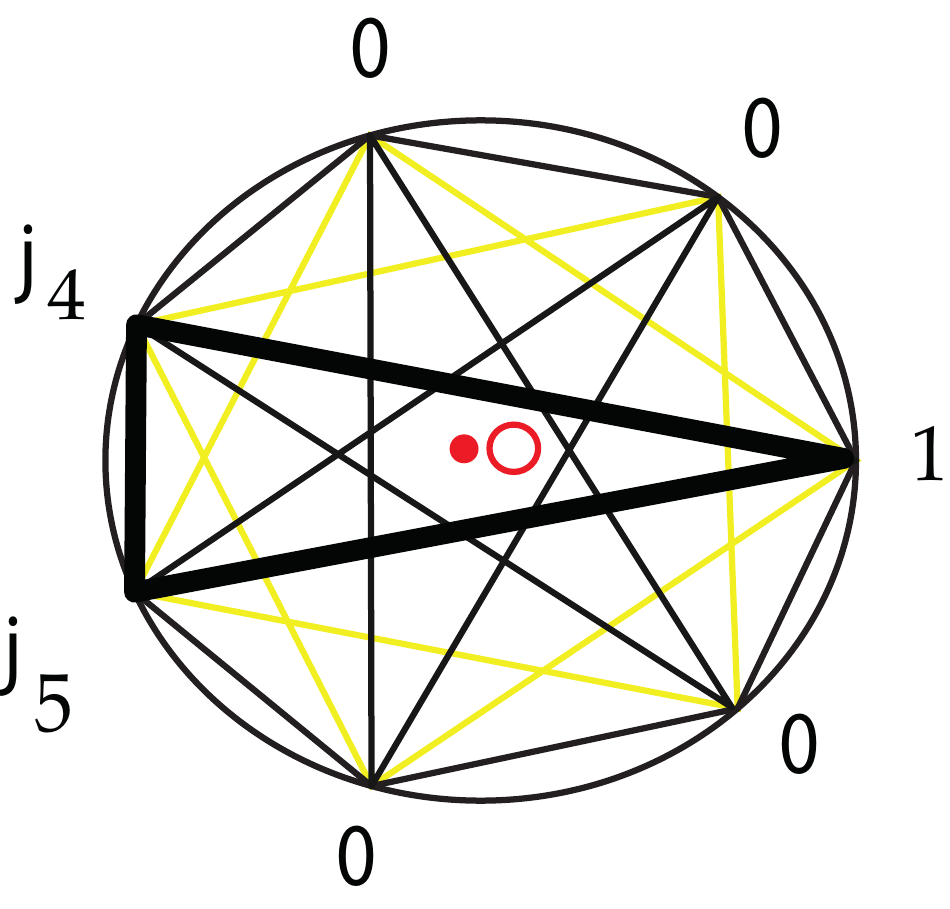}}
\end{figure}
\end{center}   
\vspace{-0.25in}
\noi The corresponding spheres are obtained by reflecting in the hyperplanes co\-rres\-pon\-ding to elements in $D_2$ and $D_{\ell+2}$, respectively. Since these sets are disjoint, the only point of intersection is the point $v$.\\

\noi Now, a neighborhood of the vertex $v$ in $\p$ looks like the first orthant of $\R^{n-3}$ where the faces $F_{L_2}$ and $F_{L_{\ell+2}}$ correspond to complementary subspaces. When reflected in all the coordinates hyperplanes of $\R^{n-3}$, one obtains a neighborhood of $v$ in $Z_+$ where those subspaces produce neighborhoods   of the two spheres. Therefore the spheres intersect transversely in that point.\\

\noi A regular neighborhood of the union of those spheres is diffeomorphic to their product minus a disk:
$$
\s^{d_2-1}\times\s^{d_{\ell+2}-1}\backslash \disc^{n-3}.
$$

\noi Joining its boundary with the boundary of $Z_+$ we see that $Z_+$ is the connected sum along the boundary of two manifolds:
$$Z_+= \s^{d_2-1}\times\s^{d_{\ell+2}-1}\backslash \disc^{n-3}\coprod X$$
where $\partial X =Z_0$ and $X$ is simply connected. Now, all the homology of $X$ comes from its boundary which again is $Z_{0}$ , since all those classes actually live in the homology of $Z$ and are the ones corresponding to the clases $D_i$ and $\tilde{D}_i$ that do not contain $n_{1}$. Those classes also exist in the homology of $Z_{0}$ and are given by the same generators, so this part of the homology of $Z_{0}$ embeds isomorphically onto the homology of $X$.\\
  
\noi If $\ell>2$, $Z_{0}$ is a connected sum of sphere products, so the homology classes of $X$ can be represented again by disjoint products $\s^p\times\disc^{n-p-3}$ and finally we construct the analog of the manifold with boundary $Q$ and the $h$-cobordism theorem gives:\\

\noi \emph{If $Z$ is simply connected of dimension at least $6$, and $n_1=1$, $\ell>2$ then}
$$
Z_{_+}=\coprod_{i=3}^{\ell+1}\left(\s^{d_i-1}\times\disc^{n-d_i-2}\right)\coprod\coprod_{i=\ell+3}^1\left (\disc^{d_i-1}\times\s^{n-d_i-2}\right)
$$
$$
\coprod \left(\s^{d_2-1}\times\s^{d_{\ell+2}-1}\backslash \disc^{n-3}\right).
$$
\noi When $n_1=1$ and $\ell=2$ we have the additional complication that restricting to $x_1=0$ takes us from the \emph{pentagonal} $Z_+$ to the \emph{triangular} $Z_0$ , which is not a connected sum but a product of three spheres and not all of its homology below the middle dimension is spherical.
\begin{center}
\begin{figure}[h]
\vspace{-0.25in}\centerline{\includegraphics[height=3.5cm]{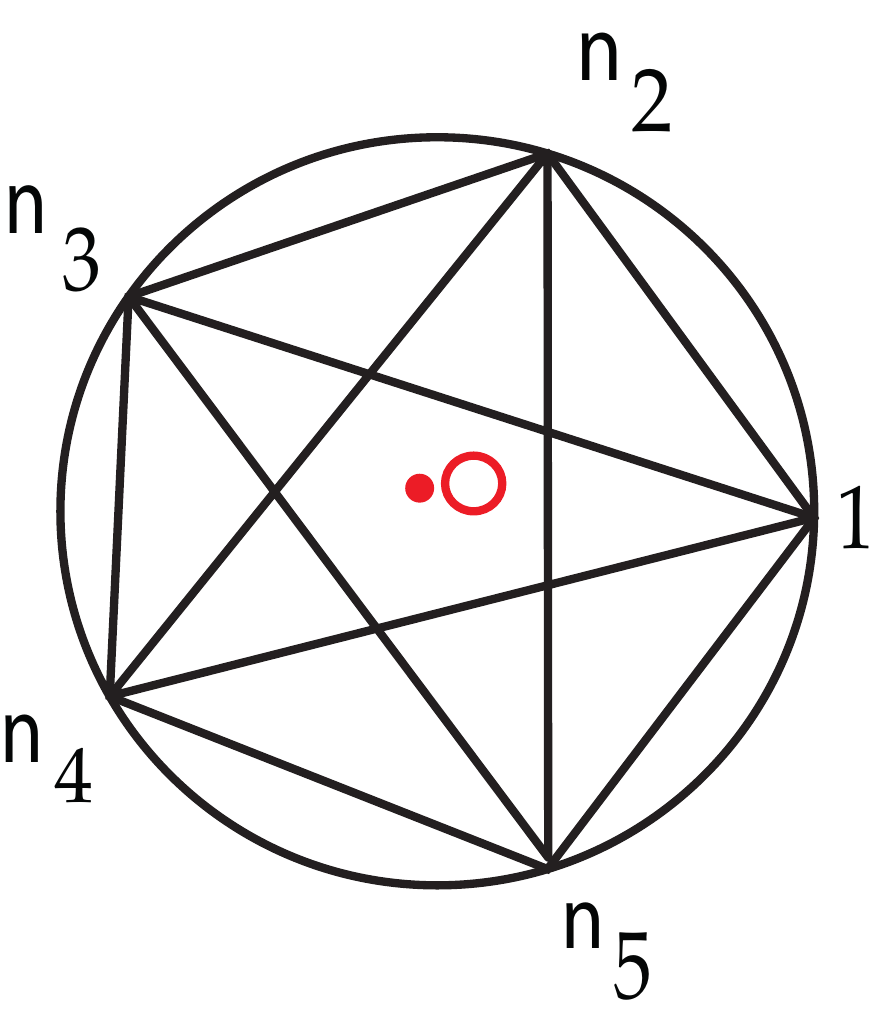}}
\end{figure}
\end{center}   
\vspace{-0.25in}
\noi The homology classes of $Z_{+}$ are those corresponding to $D_{2}$, $D_{4}$ (not coming from the boundary) and to $D_{3}$, $\tilde{D}_{1}$, $\tilde{D}_{5}$. Clearly the last ones come from the classes $[\s^{n_{3}+n{4}-1}]$, $[\s^{n_{2}-1}\times \s^{n_{3}+n{4}-1}]$ and $[\s^{n_{5}-1}\times \s^{n_{3}+n{4}-1}]$ in the boundary. This means that $X$ satisfies the hypotheses of the Lemma in section I.3 with $p=n_{2}-1$, $q=n_{5}-1$ and $m=n-3$, so we can conclude that $X$ is diffeomorphic to $\te_{n_{2}-1,n_{5}-1}^{n-3}$. We have proved all the cases of the
\begin{theo}
Let $k=2$, and consider the manifold $Z$ corresponding to the cyclic decomposition $n=n_1+\dots+n_{2\ell+1}$ and the half manifold $Z_{_+}=Z\cap\{x_1\geqslant 0\}$. When $\ell>1$ assume $Z$ and $Z_{_0}=Z\cap\{x_1=0\}$ are simply connected and the dimension of $Z$ is at least $6$. Then $Z_+$ diffeomorphic to:
\begin{itemize}
\item [a)] If $\ell=1$, the product 
$$
\s^{n_2-1}\times\s^{n_3-1}\times\disc^{n_1-1}.
$$
\item [b)] If $\ell>1$ and $n_1>1$, the connected sum along the boundary of $2\ell+1$ manifolds:
$$
\coprod_{i=2}^{\ell+2}\left(\s^{d_i-1}\times\disc^{n-d_i-2}\right)\coprod
\coprod_{i=\ell+3}^1\left(\disc^{d_i-1}\times\s^{n-d_i-2}\right).
$$
\item [c)] If $n_1=1$ and $\ell>2$, the connected sum along the boundary of $2\ell$ manifolds:
$$
\coprod_{i=3}^{\ell+1}\left(\s^{d_i-1}\times\disc^{n-d_i-2}\right)\coprod\coprod_{i=\ell+3}^1\left(\disc^{d_i-1}\times\s^{n-d_i-2}\right)
$$
$$
\coprod\left(\s^{d_2-1}\times\s^{d_{\ell+2}-1}\backslash \disc^{n-3}\right).
$$
\item [d)] If $n_1=1$ and $\ell=2$, the connected sum along the boundary of two manifolds:
$$
\left(\s^{d_2-1}\times\s^{d_4-1}\backslash \disc^{n-3}\right) \coprod \te_{n_{2}-1,n_{5}-1}^{n-3}.
$$ 
\end{itemize}
\end{theo}

\noi Theorem 3 immediately describes, under the same hypotheses, the topology of the page of the open book decomposition of $Z'$ given by Theorem 1, since this page is precisely the interior of $Z_{+}$.\\

\noi Theorem 2 about the page of the open book decomposition of the moment-angle manifold $Z^\C$ follows also, since this page is $Z_{+}$ for $Z$ the (real) intersection of quadrics corresponding to the partition $2n-1=(2n_{1}-1) + (2n_{2}) + \dots + (2n_{2\ell+1})$. In this case all the extra hypotheses of Theorem 4 hold automatically.\\

\noi Theorem 3 applies also to the topological description of some \textit{smoothings} of the cones on our intersections of quadrics (improving \cite{SLM2}). In this case the normal form is not sufficient to describe all possibilities as it was in (\cite{SLM}) where actually only the sums $d_{i}$ were needed to describe the topology or in the present work where additional information about $n_1$ only is required.\\
\newpage

\section*{Part II. Some contact structures on moment-angle mani\-folds }\label{S3}

\noi The even dimensional moment-angle manifolds and the $\lv$-manifolds have the characteristic that, except for a few, well-determined cases, do not admit symplectic structures.
We will show that the odd-dimensional moment-angle manifolds (and large families of intersections of quadrics) admit contact structures.

 \begin{theo} If $k$ is even, $\ZLam$ is a contact manifold.
  \end{theo}
 
 \noi First we show that $\ZLam$ is an almost-contact manifold. Recall that a $(2n+1)$-dimensional manifold $\vd$ is called \emph{almost contact} if its tangent bundle admits a reduction to $\su(n)\times \R$. This is seen as follows: consider the fibration $\pi: \ZLam\to \vdu(\Lam)$ with fibre the circle, given by taking the quotient by the diagonal action. Since $\vdu(\Lam)$ is a complex manifold, the foliation defined by the diagonal circle action is transversally holomorphic. Therefore,  $\ZLam$ has an atlas modeled on $\C^{n-2}\times\R$ with changes of coordinates of the charts of the form
 $$
 \left(\left(z_1,\cdots,z_{n-2}\right), t\right)  \mapsto \left(h\left(z_1,\cdots,z_{n-2},t\right), g\left(z_1,\cdots,z_{n-2},t\right)\right),
 $$ 
 \noi where $h:U\to\C^{n-2}$  and $g:U\to\R$ where $U$ is an open set in $\C^{n-2}\times\R$ and, for each
 fixed $t$ the function $\left(z_1,\cdots,z_{n-2}\right)\mapsto h\left(z_1,\cdots,z_{n-2},t\right)$ is a biholomorphism onto an open set of $\C^{n-2}\times\{t\}$. This means that the differential, in the given coordinates, is represented by a matrix of the form
  $$
 \left[\begin{array}{ccc|c}
  &  &  &  \\
   & A &  & \ast \\
  &  &  & \\ \hline
  0& \dots & 0 & r
  \end{array}\right]
$$
 \noi where $\ast$ denotes a column $(n-2)$-real vector and $A\in{\gl(n-2,\C)}$. 
 The set of matrices of the above type form a subgroup
 of $\gl(2n-3,\R)$. By Gram-Schmidt this group retracts onto $\su(n-2)\times\R$.\\
 
 \noi Now it follows from \cite{BEM} that $\ZLam$ is a contact manifold and the Theorem is proved.\\
 
 \noi  In \cite{BV} the first and third authors give a different construction, in some sense more explicit, of contact structures, not on moment-angle manifolds but on certain non-diagonal generalizations of moment-angle manifolds of the type that has been studied by  G\'omez Guti\'errez and the second author in \cite{GL}. It consists in the construction of a positive confoliation which is conductive and then uses the heat flow method described in \cite{AltWu}.\\

 \noi In a previous version of this article we had shown that the Theorem was true for many infinite families of odd-dimensional moment-angle manifolds, so we had conjectured that it would be true for all. The argument used there applies however for many other intersections of quadrics that are not moment angle manifolds, for which the proof of the previous Theorem need not apply:

\begin{theo}
There are infinitely many infinite families of odd-dimensional generic intersections of quadrics that admit contact structures.
\end{theo}

\noi First consider the odd-dimensional intersections of quadrics that are connected sums of spheres products: \\

\noi An odd dimensional product $\s^m\times\s^n$ of two spheres admits a contact structure by the following argument: let $n$ even and $m$ odd, and $n,m>2$. Without loss of generality, we suppose that $m>n$ (the other case is analogous) then $\s^m$ is an open book with binding $\s^{m-2}$ and page $\R^{m-1}$. Hence $\s^n\times\s^m$ is an open book with binding $\s^{m-2}\times\s^n$ and page $\R^{m-1}\times\s^{n}$. This page is paralellizable since
$\R\times\s^n$ already is so. Then, since $m+n-1$ is even the page has an almost complex structure. Furthermore, by hypo\-thesis, $2n\leq{n+m}$ hence by a theorem of Eliashberg \cite{Elias} the page is Stein and is the interior of an compact manifold with contact boundary $\s^{m-2}\times\s^n$. Hence by  a theorem of E. Giroux \cite{Giroux} $\s^n\times\s^m$ is a contact manifold. \\

\noi Now, it was shown by C. Meckert \cite{Meckert} and more generally by Weinstein \cite{Weinstein} (see also \cite{Elias}) that the connected sum of contact manifolds of the same dimension is a contact manifold. Therefore all odd dimensional connected sums of sphere products admit contact structures.\\

\noi Additionally, it was proved by F. Bourgeois in \cite{Bour} (see also Theorem 10 in \cite{Giroux}) that if a closed manifold $\vd$ admits a contact structure, then so does $\vd\times \T$. Therefore, all moment-angle manifolds of the form $Z\times\T^{2m}$, where $Z$ is a connected sum of sphere products, admit contact structures.\\

\noi For every case where $Z$ is a connected sum of sphere products we have an infinite family obtained by applying construction $Z\mapsto Z'$ an infinite number of times and in the different coordinates (as well as other operations). The basic cases from which to start these infinite families constitute also a large set (see \cite {LoGli}, remarks after Theorem 2.4) and the estimates on their number in each dimension keep growing. Adding to those varieties their products with tori we obtain an even larger set of cases where an odd-dimensional $Z$ admits a contact structure.\\

\noi Another interesting fact is that most of them (including moment-angle mani\-folds) also have an open book decomposition. However, for these open book decompositions there does not exist a contact form which is supported in the open book decomposition  like in Giroux's theorem because the pages are not Weinstein manifolds (i.e manifolds of dimension $2n$ with a Morse function with indices of critical points lesser or equal to $n$). It is possible however that the pages of the book decomposition admit Liouville structures in which case one could apply the techniques of D. McDuff (\cite{McDuff}) and P. Seidel (\cite{Seidel}) to obtain contact structures.

\noindent Yadira Barreto. {\tt Universidad Aeron\'autica en Quer\'etaro}.
Carretera Estatal Quer\'etaro-Tequisquiapan 22154, 76270, Col\'on, Qro.

\noindent {\it E-mail address:} yadira.barreto@unaq.edu.mx  

\vskip.3cm
\noindent Santiago L\'opez de Medrano . {\tt Instituto de Matem\'aticas}. Universidad Nacional Au\-t\'o\-no\-ma de M\'exico.
Ciudad Universitaria. Coyoac\'an 04510, M\'exico, D. F. M\'exico. 

\noindent {\it E-mail address:} santiago@matem.unam.mx

\vskip .3cm

\noindent Alberto Verjovsky. {\tt Instituto de Matem\'aticas, Unidad Cuernavaca}. Universidad Nacional Au\-t\'o\-no\-ma de M\'exico.
Av. Universidad s/n, Col. Lomas de Chamilpa. Cuernavaca, Morelos, M\'exico, 62209.

\noindent {\it E-mail address:} alberto@matcuer.unam.mx


\begin{thebibliography}{100}

\bibitem{AltWu} Altschuler, S. J. and Wu, L. F., ``On deforming confoliations", \\
\emph{J. Differential Geometry}, Vol. 54, pp. 75-97, 2000.

\bibitem{BBCG} Bahri, A., Bendersky, M., Cohen, F. R. and Gitler, S.,  ``The polyhedral product functor: a method of decomposition for moment-angle complexes, arrangements and related spaces"., \emph{Adv. Math.}, Vol. 225, no. 3, pp. 1634-1668, 2010.

\bibitem{BLV} Barreto, Y., L\'opez de Medrano, S. and  Verjovsky, A., ``Open Book Structures on Moment-Angle Manifolds $Z^{\mathbb C}(\Lambda)$ and Higher Dimensional Contact Manifolds". arXiv:1303.2671.

\bibitem{BV} Barreto, Y. and  Verjovsky, A.,  ``Moment-angle manifolds, intersection of quadrics and higher dimensional contact manifolds", \emph{Moscow Mathematical Journal.}, Vol. 14, No. 4, pp. 669-696, 2014.

\bibitem{Ba}  Baskakov, I. V., ``Massey triple products in the cohomology of moment-angle complexes", Russian Math. Surveys 58 (2003), 1039-1041.

 \bibitem{BEM} Borman, M.S., Eliashberg, Y. and Murphy, E.,``Existence and classification of overtwisted contact structures in all dimensions", \emph{Acta Mathematica}, 215 (2015), 281-361.
 
\bibitem{BM} Bosio, F and Meersseman, L., ``Real quadrics in $\C^n$, complex manifolds and convex polytopes", \emph{Acta Math.}, Vol. 197, pp.  53-127, 2006.

\bibitem{Bour} Bourgeois, F., ``Odd dimensional tori are contact manifolds", \emph{International Mathematics Research Notices}, no. 30, pp. 1571-1574, 2002.
 
\bibitem{Elias} Eliashberg, Y., ``Topological characterization of Stein manifolds of dimension $>2$ ," \emph{International Journal of  Mathematical}, Vol. 1, No. 1, pp. 29-46,  1990.

\bibitem{Giroux} Giroux, E., ``Geometrie de contact: de la dimension trois vers les dimensions superieures", \emph{ICM}, Vol. II, pp. 405-414, 2002.

\bibitem{GiMo} Giroux, E. and Mohsen, J. P., ``Structures de contact et fibrations symplectiques sur le cercle", \emph {in  process}.

\bibitem {LoGli} Gitler, S. and L\'opez de Medrano, S., ``Intersections of quadrics, moment-angle manifolds and connected sums", \emph{Geom. Topol.}  Vol. 17, no. 3, pp. 1497-1534, 2013.

\bibitem{GL} G\'omez Guti\'errez, V. and L\'opez de Medrano, S., ``The topology of the intersections of quadrics II", \emph{Bolet\'in de la Sociedad Matem\'atica Mexicana}, Vol. 20, pp. 237-255, 2014. 

\bibitem{Gray} Gray, J. W., ``Some global properties of contact structures", \emph{Annals of Mathematics}, Vol. 69, no. 2, pp. 421-450, 1959.

\bibitem{H} Hirsch, M.W., ``Smooth Regular Neighborhoods", \emph{Annals of Mathematics}, Vol. 76, No. 3 (Nov., 1962), pp. 524-530.

\bibitem{SLM} L\'{o}pez de Medrano, S., ``The topology of the intersection of quadrics in $\R^{^n}$", \emph{Springer-Verlag. Lectures Notes in Mathematics}, Vol. 1370,  pp. 280-292, 1989.

\bibitem{SLM2} L\'{o}pez de Medrano, S., ``Singularities of homogeneous quadratic mappings", \emph{Revista de la Real Academia de Ciencias Exactas, Físicas y Naturales.}, Serie A, Matem\'aticas, Vol. 108 (2014),  pp. 95-112, Springer-Verlag.

\bibitem{LV} L\'{o}pez de Medrano, S. and  Verjovsky, A., ``A new family of complex, compact, nonsymplectic manifolds", Bol. Soc. Brasil. Mat., 28 (1997), 253-269.

\bibitem{LM} Lutz, R. and Meckert, C., ``Structures de contact sur certaines sph\`eres exotiques", \emph{C. R. Acad. Sci. Paris S\'er. A-B }, Vol. 282, pp. A591-A593, 1976.

\bibitem{MMP} Mart\'inez, D., Mu\~noz, V. and Presas, F., ``Open book decompositions for almost contact manifolds", \emph{Proceedings of the XI Fall Workshop on Geometry and Physics, Publicaciones de la RSME}, Vol. 6, pp. 131-149, 2004.

\bibitem{McDuff} McDuff, D., ``Symplectic manifolds with contact type boundaries", \emph{Invent. Math.}, Vol. 103, no. 3, pp. 651-671, 1991.

\bibitem{Meckert} Meckert, C., ``Forme de contact sur la somme connexe de deux vari\'et\'es de contact de dimension impare", \emph{Annales De L'Institut Fourier }, Vol. 32, no. 3 pp. 251-260, 1982.

\bibitem{Meer} Meersseman, L., ``A new geometric construction of compact complex manifolds in any dimension", Math. Ann., 317 (2000), 79-115.

\bibitem{Seidel}  Seidel, P., ``Simple examples of distinct Liouville type symplectic structures",  \emph{J. Topol. Anal. }, no. 1, pp,  1-5,  2011.

\bibitem{Weinstein} Weinstein, A., ``Contact surgery and symplectic handlebodies", \emph{Hokkaido Math. J.}, Vol. 620, no. 2,  pp. 241-251, 1991.

\end{thebibliography}
\end{document}